\documentclass[numbers,final]{my_aipproc}

\layoutstyle{8x11single}

\begin{document}

\title{Insecticide control 
in a Dengue epidemics model\thanks{Accepted 28/07/2010 
to Proc. 8th Int. Conf. Numerical Analysis and Applied Mathematics
(ICNAAM 2010), Rhodes, Greece, 19-25 Sept 2010.}}



\author{Helena Sofia Rodrigues}{
  address={School of Business Studies, Viana do Castelo Polytechnic Institute, Portugal\\
  sofiarodrigues@esce.ipvc.pt}
}

\author{M. Teresa T. Monteiro}{
  address={Department of Production and Systems, University of Minho, Portugal\\
    tm@dps.uminho.pt}
}

\author{Delfim F. M. Torres}{
  address={Department of Mathematics, University of Aveiro, Portugal\\
    delfim@ua.pt}
}


\begin{abstract}
A model for the transmission of dengue disease is presented. 
It consists of eight mutually-exclusive compartments representing 
the human and vector dynamics. It also includes a control parameter 
(insecticide) in order to fight the mosquitoes. 
The main goal of this work is to investigate the best way 
to apply the control in order to effectively reduce the number 
of infected humans and mosquitoes. A case study, using data 
of the outbreak that occurred in 2009 in Cape Verde, is presented.
\end{abstract}

\maketitle


\section{Introduction}

Dengue is a mosquito-borne infection mostly found in tropical and sub-tropical 
climates worldwide, mostly in urban and semi-urban areas. It can cause a severe 
flu-like illness, and sometimes a potentially lethal complication called 
dengue haemorrhagic fever. About 40\% of the world's population are now at risk.
The spread of dengue is attributed to expanding geographic distribution of the four 
dengue viruses and their mosquito vectors, the most important of which is the predominantly 
urban species \emph{Aedes aegypti}. The life cycle of a mosquito presents four distinct stages: 
egg, larva, pupa and adult. In the case of \emph{Aedes aegypti} the first three stages 
take place in or near water while air is the medium for the adult stage. 
The adult stage of the mosquito, in the urban environment, 
is considered to last an average of eleven days.

The paper is organized as follows. The next section presents a mathematical model 
of the interaction between human and mosquito populations. Then, the numerical 
experiments using different strategies for the insecticide administration are reported. 
In the last section the conclusions are presented.


\section{The mathematical model}

Considering the work of \cite{Rodrigues2009}, the relationship between humans
and mosquitoes are now rather complex, taking into account the model presented 
in \cite{Dumont2008}. The novelty in this paper is the presence of the control
parameter related to adult mosquito insecticide \cite{MathCM2010}.

The notation used in our mathematical model includes 
four epidemiological states for humans (index $h$):

\begin{quote}
\begin{tabular}{ll}
$S_h(t)$ & susceptible (individuals who can contract
the disease)\\
$E_h(t)$ & exposed (individuals who have been infected by the parasite\\
& but are not yet able to transmit to others)\\
$I_h(t)$ & infected (individuals capable of transmitting the disease to others)\\
$R_h(t)$ & resistant (individuals who have acquired immunity)\\
\end{tabular}
\end{quote}

It is assumed that the total human population $(N_h)$ is constant,
so, $N_h=S_h+E_h+I_h+R_h$. There are also other four state variables 
related to the female mosquitoes (index $m$) 
(the male mosquitoes are not considered in this study
because they do not bite humans and consequently they do not
influence the dynamics of the disease):

\begin{quote}
\begin{tabular}{ll}
$A_m(t)$& aquatic phase (that includes the egg, larva and pupa stages)\\
$S_m(t)$& susceptible (mosquitoes that are able to contract the disease)\\
$E_m(t)$& exposed (mosquitoes that are infected but are not yet able \\
& to transmit to humans)\\
$I_m(t)$& infected (mosquitoes capable of transmitting the disease to humans)\\
\end{tabular}
\end{quote}

In order to analyze the effects of campaigns to combat the mosquito, 
there is also a control variable:
\begin{quote}
\begin{tabular}{ll}
$c(t)$& level of insecticide campaigns\\
\end{tabular}
\end{quote}

Some assumptions are made in this model: the total human population ($N_h$) is constant, 
which means that we do not consider births and deaths; there is no immigration 
of infected individuals to the human population; the population is homogeneous; 
the coefficient of transmission of the disease is fixed and do not vary seasonally; 
both human and mosquitoes are assumed to be born susceptible; 
there is no natural protection; for the mosquito there is no resistant phase, 
due to its short lifetime.

The Dengue epidemic can be modelled by the following nonlinear
time-varying state equations for Human Population and vector population, respectively:
\begin{equation}\label{general_model}
\begin{tabular}{l}
$
\left\{
\begin{array}{l}
\frac{dS_h}{dt}(t) = \mu_h N_h - (B\beta_{mh}\frac{I_m}{N_h}+\mu_h) S_h\\
\frac{dE_h}{dt}(t) = B\beta_{mh}\frac{Im}{N_h}S_h - (\nu_h + \mu_h )E_h\\
\frac{dI_h}{dt}(t) = \nu_h E_h -(\eta_h  +\mu_h) I_h\\
\frac{dR_h}{dt}(t) = \eta_h I_h - \mu_h R_h\\
\end{array}
\right. $\\
\end{tabular}
\begin{tabular}{l}
$
\left\{
\begin{array}{l}
\frac{dA_m}{dt}(t) = \mu_b (1-\frac{A_m}{K})(S_m+E_m+I_m)-(\eta_A+\mu_A) A_m\\
\frac{dS_m}{dt}(t) = -(B \beta_{hm}\frac{I_h}{N_h}+\mu_m) S_m+\eta_A A_m-c S_m\\
\frac{dE_m}{dt}(t) = B \beta_{hm}\frac{I_h}{N_h}S_m-(\mu_m + \eta_m) E_m-c E_m\\
\frac{dI_m}{dt}(t) = \eta_m E_m -\mu_m I_m - c I_m\\
\end{array}
\right. $\\
\end{tabular}
\end{equation}
with the initial conditions: $S_h(0)=S_{h0},\;  E_h(0)=E_{h0},\; I_h(0)=I_{h0},\;
R_h(0)=R_{h0}, \; A_m(0)=A_{m0}, \; S_{m}(0)=S_{m0},\;
E_m(0)=E_{m0},\;  I_m(0)=I_{m0}.$

Notice that the equation related to the aquatic phase for the mosquito does not have
the control variable $c$, because this kind of insecticide does not produce
effects in this stage.

The parameters used in the model are:
\begin{quote}
\begin{tabular}{llll}
$N_h$ & total human population & $\mu_{A}$ & natural mortality of larvae (per day)\\
$B$ & average daily biting (per day)& $\eta_{A}$ & maturation rate from larvae to adult (per day)\\
$\beta_{mh}$ & transmission probability from $I_m$ (per bite) &$1/\eta_{m}$ & extrinsic incubation period (in days) \\
$\beta_{hm}$ & transmission probability from $I_h$ (per bite) & $1/\nu_{h}$ & intrinsic incubation period (in days)\\
$1/\mu_{h}$ & average lifespan of humans (in days) & $m$ & female mosquitoes per human\\
$1/\eta_{h}$ & mean viremic period (in days)& $k$ & number of larvae per human\\
$1/\mu_{m}$ & average lifespan of adult mosquitoes (in days) & $K$ & maximal capacity of larvae\\
$\mu_{b}$ & number of eggs at each deposit per capita (per day) & & \\
\end{tabular}
\end{quote}

\begin{figure}[ptbh]
  \includegraphics [scale=0.25]{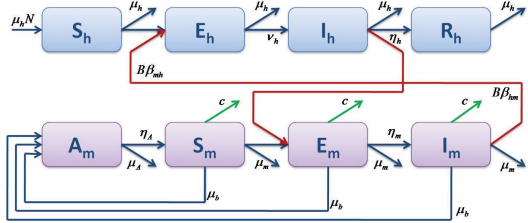}\\
   {\caption{\label{model}  Interaction between human and mosquito}}
\end{figure}

The Figure~\ref{model} shows the relation between human 
and mosquito and the corresponding model parameters.


\section{Numerical experiments}

The values related to humans originate from data collected 
in an infected period in Cape Verde \cite{CDC2010}. However,
as it was the first outbreak that happened in the archipelago, 
it was not possible to collect any data for the mosquito. 
Thus, for the \emph{Aedes aegypti} we have selected 
information from Brazil where dengue is already 
a long known reality \cite{Thome2010,Yang2009}.

The numerical tests were carried out using Scilab \cite{Scilab} with
the following values \cite{Thome2010,Yang2009}:
$N_h=480000$, $B=1$, $\beta_{mh}=0.375$,
$\beta_{hm}=0.375$, $\mu_{h}=1/(71\times 365)$, $\eta_{h}=1/3$,
$\mu_{m}=1/11$, $\mu_{b}=6$, $\mu_{A}=1/4$, $\eta_A=0.08$,
$\eta_m=1/11$, $\nu_h=1/4$, $m=6$,
$k=3$, $K=k\times N_h$.
The initial conditions for the problem were: $S_{h0}=N_h-E_{h0}-I_{h0}$,
$E_{h0}=216$, $I_{h0}=434$, $R_{h0}=0$, $A_{m0}=k\times N_h$,
$S_{m0}=m\times N_h$, $E_{m0}=0$, $I_{m0}=0$. The final time was
$t_f=84$ days.
\begin{figure}[ptbh]
\centering
\begin{minipage}[t]{0.48\linewidth}
\centering
\includegraphics[scale=0.52]{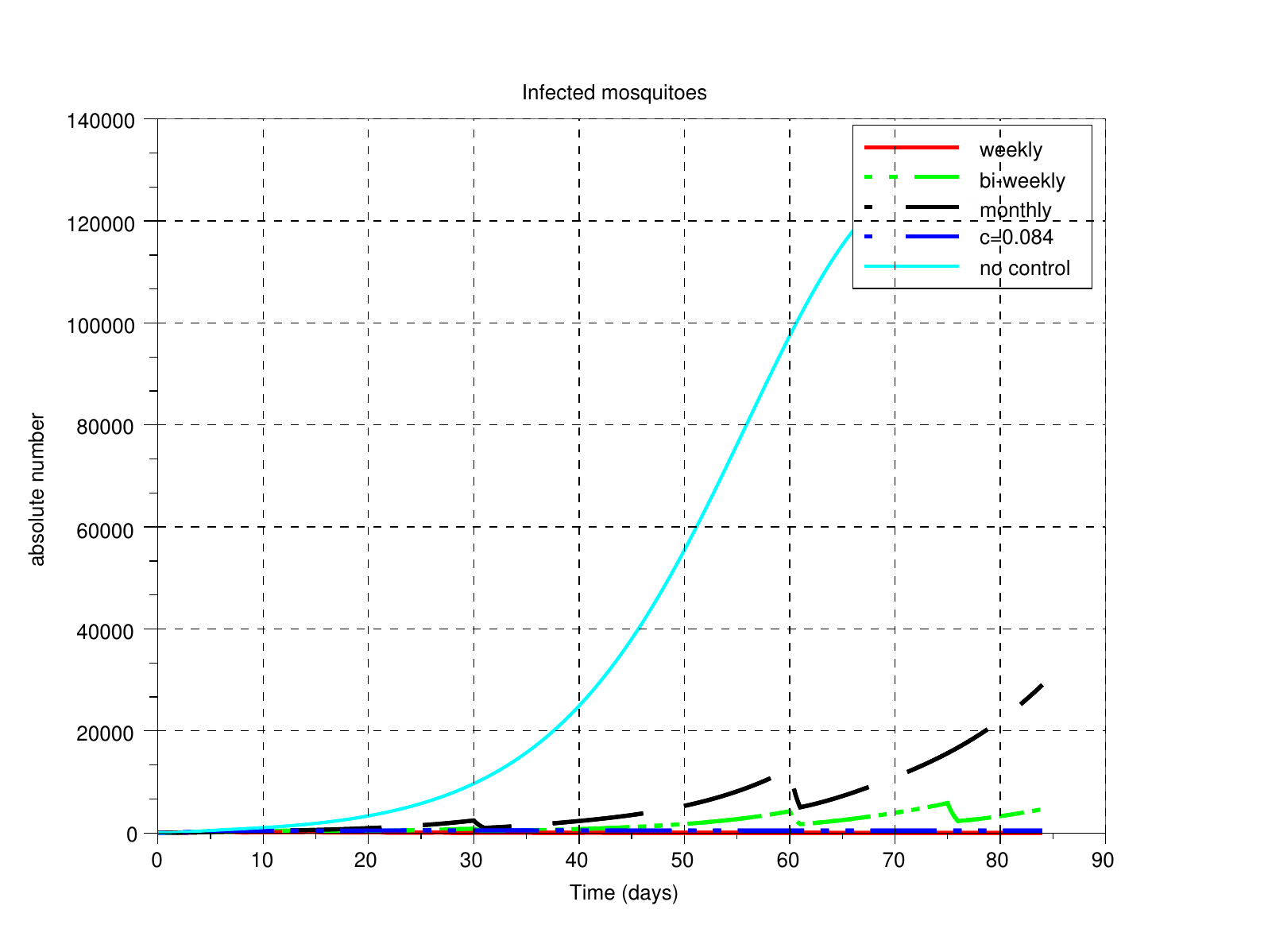}
\end{minipage}\hspace*{\fill}
\begin{minipage}[t]{0.48\linewidth}
\centering
\includegraphics[scale=0.52]{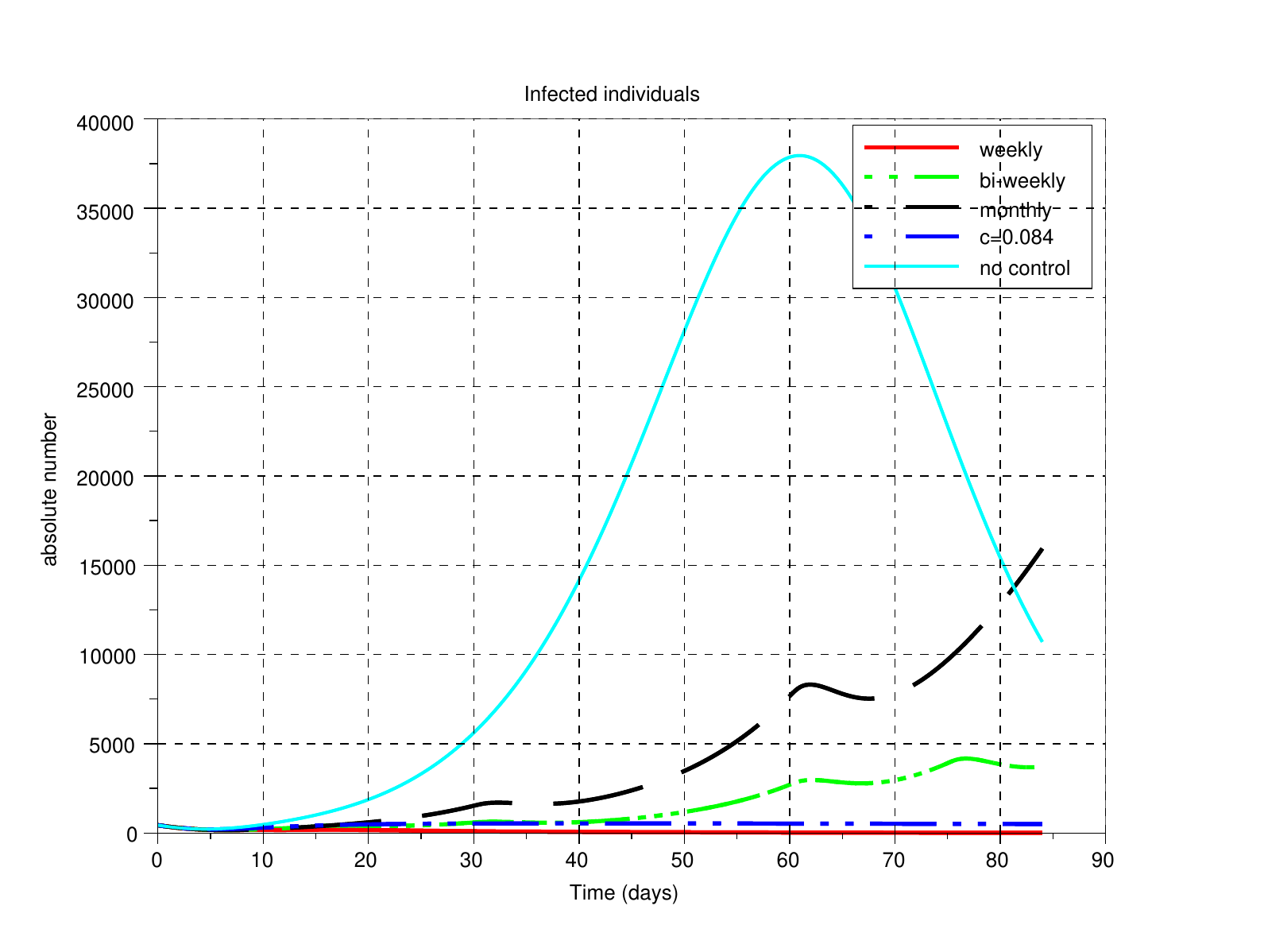}
\end{minipage}
\caption{\label{control} Infected mosquitoes and individuals}
\end{figure}

In the literature it has been proved that a disease-free equilibria (DFE) 
is locally asymptotically stable, whenever a certain epidemiological threshold, 
known as the \emph{basic reproduction number}, is less than one. 
The \emph{basic reproduction number} of the disease represents the expected 
number of secondary cases produced in a completed susceptible population, 
by a typical infected individual during its entire period of infectiousness 
\cite{Hethcote2000}. In a recent work \cite{CMMSE2010} it was proved that
if a constant minimum level of insecticide is applied ($c=0.084$), 
it is possible to maintain the basic reproduction number below unity, guaranteing the DFE.

\begin{figure}[ptbh]
\centering
\begin{minipage}[t]{0.48\linewidth}
\centering
\includegraphics[scale=0.52]{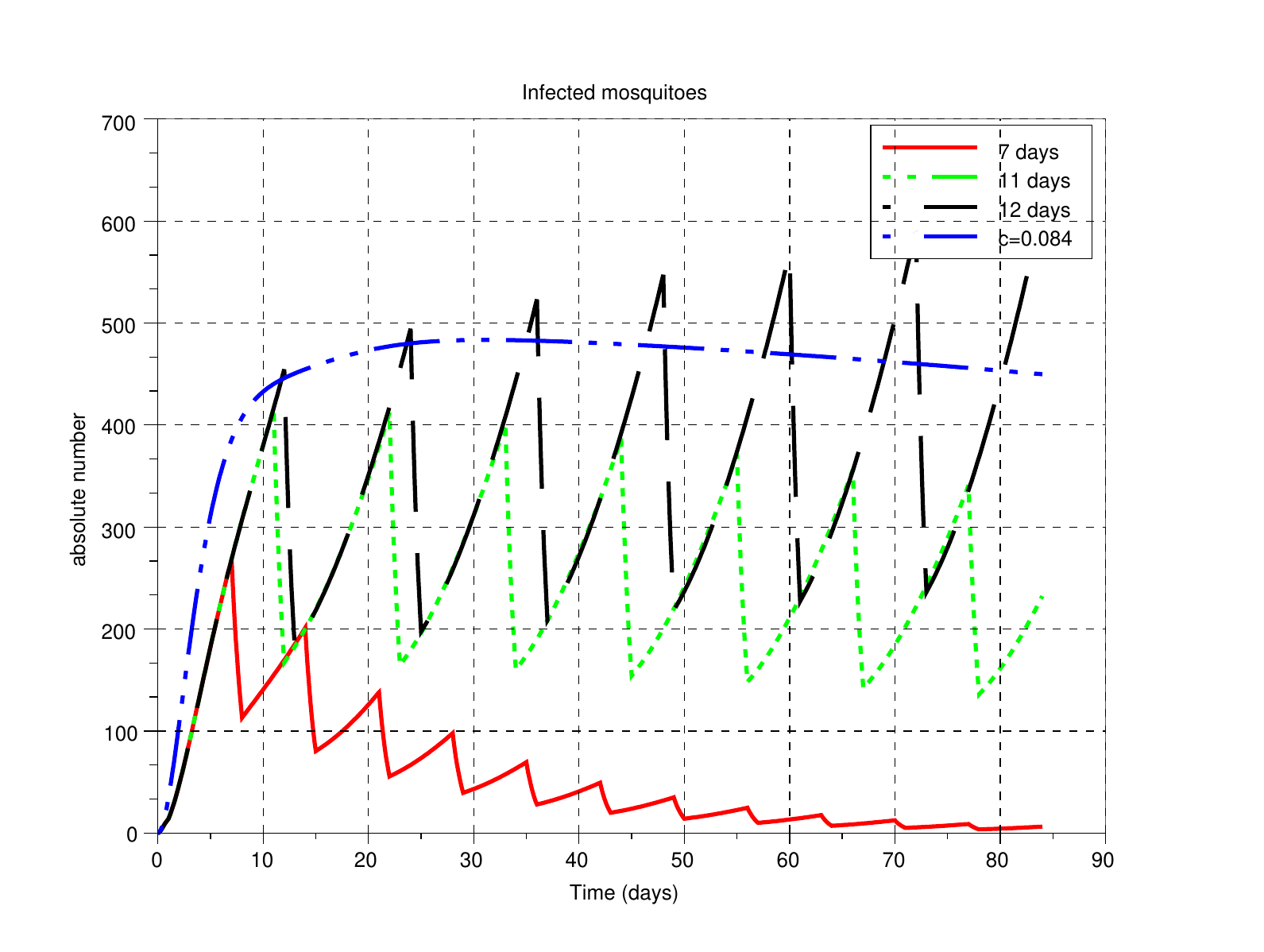}
\end{minipage}\hspace*{\fill}
\begin{minipage}[t]{0.48\linewidth}
\centering
\includegraphics[scale=0.52]{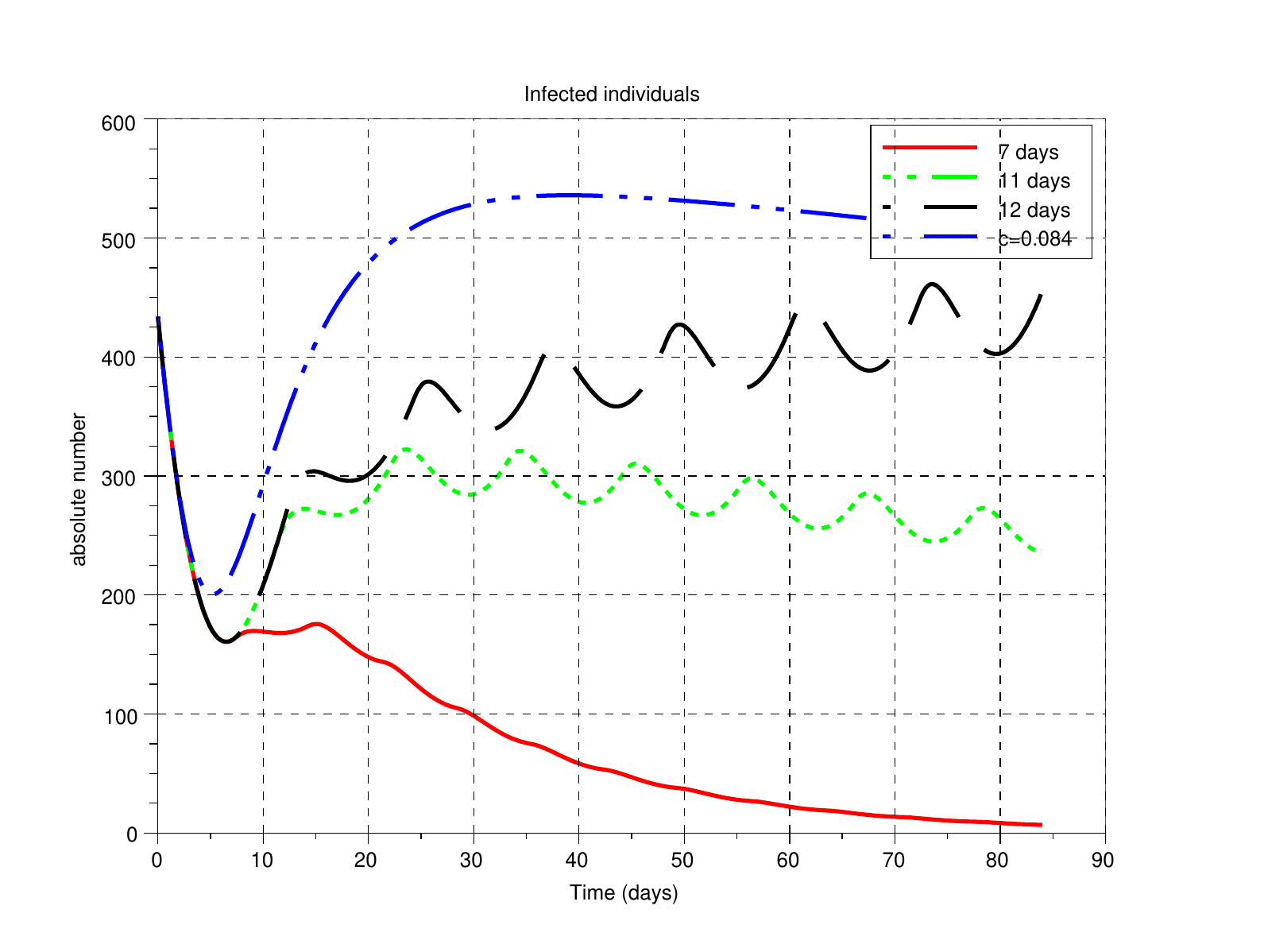}
{\caption{\label{periodicity} Infected mosquitoes and individuals \emph{vs} periodicity}}
\end{minipage}
\end{figure}
To solve the system (\ref{general_model}), in a first step,  several strategies 
of control application were used: three different frequencies 
(weekly, bi-weekly and monthly) for control application, constant control 
($c=0.084$ in \cite{CMMSE2010}) and no control ($c=0$). The weekly/bi-weekly/monthly
frequency means that during one day per week/bi-week/month the whole ($100\%$) 
capacity of insecticide ($c=1$) is used during all day. 
Constant control strategy ($c=0.084$) consists in the application of $8.4\%$  
capacity of insecticide 24 hours per day all the time (84 days). In this work, 
the amount of insecticide is an adimensional value and must be considered in relative terms.

Figure~\ref{control} shows the results of these strategies regarding infected mosquitoes 
and individuals. Without control, the number of infected mosquitoes 
and individuals increases expressively.

Realizing the influence of the insecticide control, further tests were carried out 
to find the optimum periodicity of administration which, from gathered results, 
must rest between one and two weeks. The second phase of numerical tests, 
Figure~\ref{periodicity}, considers four situations: 7 days, 11 days, 12 days 
and continuously $c=0.084$. To guarantee the DFE, the curves must remain below 
the one corresponding to $c=0.084$.

The amount of insecticide, and when to apply it, are important factors for outbreak control. 
Table~\ref{result} reports the total amount of insecticide 
used in each version during the 84 days.

\begin{table}[ptbh]
\begin{tabular}{ccccccc}
\hline
& 7 days & 11 days & 12 days & 15 days & 30 days & $c=0.084$ \\
\hline
\hline
insecticide amount  & 12& 8 & 7 & 6 & 3 & 7.056  \\\hline
\end{tabular}
\caption{Insecticide cost}
\label{result}
\end{table}


\section{Conclusions}

The numerical tests conclude that the best strategy for the infected reduction 
is the weekly administration, however it is also the most expensive one (insecticide cost). 
The best result obtained is between 11 and 12 days, with the insecticide 
amount in the closed interval from 7 to 8, confirming the amount for $c=0.084$ 
in \cite{CMMSE2010}. The 11 or 12 days between applications can be directly 
related to the span of adult stage for the mosquitoes, 
an average of eleven days in an urban environment.

In this work the insecticide administration time was considered continuous 
in each day (24 hours per day). As a future work, an optimization problem 
will be formulated and solved in order to find the best plan in terms 
of periodicity and duration (limited number of hours per day). 
The plan must consider the advantages to apply the insecticide only during the night.


\begin{theacknowledgments}
Partially supported by the Portuguese Foundation 
for Science and Technology (FCT) through 
the PhD Grant SFRH/BD/33384/2008 (Rodrigues).
\end{theacknowledgments}




\begin{thebibliography}{9}

\bibitem{CDC2010}
\emph{CDC 2010}: 
{http://www.cdc.gov/dengue/}, {April}, {2010}.

\bibitem{Dumont2008}
Y.~Dumont, F.~Chiroleu, and C.~Domerg, 
\emph{On a temporal model for the Chikungunya disease: 
modeling, theory and numerics}, Math. Biosci., 
\textbf{213}-1, 80--91, (2008).

\bibitem{Hethcote2000}
H.~W.~ Hethcote, 
\emph{The mathematics of infectious diseases}, 
SIAM Rev., \textbf{42}-4, 599--653, (2000).

\bibitem{Rodrigues2009}
H.~S.~Rodrigues, M.~T.~T.~ Monteiro, and D.~F.~M.~ Torres, 
\emph{Optimization of Dengue Epidemics: A Test Case with Different Discretization Schemes}, 
in \emph{Numerical Analysis and Applied Mathematics}, 
{AIP Conference Proceedings 1168}, American Institute of Physics, {2009},  pp.~1385--1388.
{\tt arXiv:1001.3303}

\bibitem{CMMSE2010}
H.~S.~ Rodrigues, M.~T.~T.~ Monteiro, D.~F.~M.~ Torres, and A.~ Zinober, 
\emph{Control of dengue disease: a case study in Cape Verde}, 
in Proceedings of the 10th International Conference on Computational 
and Mathematical Methods in Science and Engineering, 
CMMSE June 2010, Almeria, 2010, pp. 816--822. 
{\tt arXiv:1006.5931}

\bibitem{MathCM2010}
H.~S.~Rodrigues, M.~T.~T.~ Monteiro, and D.~F.~M.~ Torres, 
\emph{Dynamics of Dengue epidemics when using optimal control}, 
Mathematical and Computer Modelling, 2010, 
DOI: 10.1016/j.mcm.2010.06.034.
{\tt arXiv:1006.4392}

\bibitem{Scilab}
\emph{Scilab}: http://www.scilab.org/, {May}, {2010}.

\bibitem{Thome2010}
R.~C.~ Thom\'e, H.~M.~ Yang, and L.~ Esteva,
\emph{Optimal control of Aedes aegypti mosquitoes 
by the sterile insect technique and insecticide},
Math. Biosci. \textbf{223}-1,12--23, (2010).

\bibitem{Yang2009}
H.~M.~ Yang, M.~L.~G.~ Macoris,  K.~C.~ Galvani, M.~T.~M.~Andrighett, and D.~M.~V.~Wanderley,
\emph{Assessing the effects of temperature on dengue transmission}, 
Epidemiol. Infect., Cambridge University Press, 2009.

\end{thebibliography}
\end{document}